\documentclass[sn-mathphys]{sn-jnl}
\jyear{2022}%

\theoremstyle{thmstyleone}%
\newtheorem{theorem}{Theorem}

\def\intt{[\kern .3pt 0,1]}
\def\twon{2\kern .4pt n}
\def\qed{~\hbox{\vrule width 3pt depth 4pt height 2.5pt}}

\theoremstyle{thmstylethree}%
\def\cmax{\max_k^{} |c_k^{}|}
\def\pput(#1,#2)#3{\noindent\smash{\raise#2pt\hbox to 0pt
   {\kern #1pt #3\hss}}\ignorespaces}

\begin{document}

\title[Spectacularly large expansion coefficients]{Spectacularly
large expansion coefficients in M\"untz's theorem}

\author*{\fnm{Lloyd N.} \sur{Trefethen}}\email{trefethen@maths.ox.ac.uk}

\affil*{\orgdiv{Mathematical Institute},
\orgname{University of Oxford},
\orgaddress{\street{Woodstock Road}, \city{Oxford}, \postcode{OX2 6GG},
\country{United Kingdom}}}

\abstract{M\"untz's theorem asserts, for example, that the even powers $1, x^2, x^4,\dots$
are dense in $C(\intt)$.  We show that the associated expansions
are so inefficient as to have no conceivable relevance to any actual computation.
For example, approximating
$f(x)=x$ to accuracy $\varepsilon = 10^{-6}$ in this basis requires powers larger than $x^{280{,}000}$ and coefficients
larger than $10^{107{,}000}$.  We present a theorem establishing 
exponential growth of coefficients with respect to $1/\varepsilon$.}



\maketitle


\section{Introduction}\label{sec1}
The M\"untz approximation theorem, conjectured by Bernstein in
1912~\cite{b12} and proved by M\"untz in 1914~\cite{muentz}, is a beautiful result. 
Suppose we are interested in continuous functions on $\intt$,
i.e., $f\in C(\intt)$, and we want to approximate them as linear
combinations of the monomials
$$
x^{a_1^{}}, x^{a_2^{}}, \dots,
$$
where $\{a_k^{}\}$ is a set of exponents (not necessarily
integers) satisfying
$$
0 = a_0^{} < a_1^{} < a_2^{} < \cdots .
$$
Certainly this is possible if $a_k^{} = k$ for each $k$, by the
Weierstrass approximation theorem~\cite{atap}, but what
if the set of powers is sparser?  For example, are the even
powers
\begin{equation}
1, x^2, x^4, \dots
\label{evens}
\end{equation}
dense in $C(\intt)$?  The theorem characterizes all suitable sets of exponents:

\begin{theorem}[M\"untz approximation theorem]\label{thm1}
The family $\{x^{a_k^{}}\}$ is dense in $C([0,1])$ if and only
if 
\begin{equation}
\sum_{k=1}^\infty {1\over a_k^{}} = \infty .
\label{muentzcond}
\end{equation}
\end{theorem}

\noindent
Thus the set (\ref{evens}) easily qualifies, as do many other
collections of exponents, such as the primes:
\begin{equation}
1, x^2, x^3, \dots, x^{89}, x^{97}, \dots.
\label{primes}
\end{equation}
For discussions
of the theorem with proofs, see~\cite{almira}, ~\cite{be}, and~\cite{lz}.
After 1914, M\"untz's theorem was generalized by
Sz\'asz and others. 

As a numerical analyst, I work with algorithms based on
expanding functions in nonorthogonal bases, a powerful technique in
certain contexts~\cite{fna1,rect}.
This led me to consider M\"untz's theorem from a computational
angle, and what emerged is startling.  To make the point, it
is enough to consider a particular case of
what might be regarded as the most basic
nontrivial M\"untz approximation.  The name ``E'' alludes to the use
of even powers.

\medskip

{\em Problem E.  Given $\varepsilon>0$, find an integer $n\ge 0$ and coefficients
$c_0^{}, \dots, c_n^{}$ such that 
\begin{equation}
\Bigl|\kern 1pt x - \sum_{k=0}^n c_k^{} x^{2k}\Bigr| \le \varepsilon, \qquad
x \in [\kern .3pt 0,1].
\label{prob}
\end{equation}
}

\noindent We shall prove:

\begin{theorem}\label{thm2}
If\/ $\varepsilon < 1/2$, then any solution of Problem E has
\begin{equation}
n > {1\over 20\kern .5pt \varepsilon}
\label{Nbound}
\end{equation}
and
\begin{equation}
\cmax > 0.75 \kern .7pt \varepsilon \kern .7pt 2^{\kern .8pt 1/(40 \kern .4pt \varepsilon)} . 
\label{cbound}
\end{equation}
\end{theorem}

\noindent Actually, I believe the following sharper bounds hold:
\begin{equation}
n > {1\over 8\varepsilon},
\label{Nboundest}
\end{equation}
\begin{equation}
\cmax > {(1+\sqrt 2 \kern 1pt)^{\twon} \over 16\kern .5pt n^{1.5}}.
\label{cboundest}
\end{equation}
For accuracy $\varepsilon = 10^{-6}$, my estimate is that
one needs $n> 140{,}000$ and $\cmax > 10^{107{,}000}$.
In such an expansion, the enormous 
coefficients have oscillating signs, so that they cancel
almost exactly (namely to one part in $10^{107{,}000}$).  
On a computer in floating-point arithmetic, all information
will be lost unless one works in a precision of more than $107{,}000$ digits.
(The usual precision is 16 digits.)

\section{Proof}

Problem E is 
equivalent to the more familiar problem of approximation of $|x|$ on $[-1,1]$:
\begin{equation}
\Bigl| |x| - \sum_{k=0}^{2\kern .3pt n} c_k^{} x^{2k}\Bigr| \le \varepsilon, \qquad
 x \in [-1,1].
\label{prob2}
\end{equation}
Since Lebesgue first used approximations of $|x|$ for a proof of
the Weierstrass approximation theorem at age 23 in 1898, a great deal has 
been learned about this problem, as recounted in
chapter~25 of~\cite{atap}.  In particular, Bernstein's 1914 paper~\cite{bernstein}
was a landmark contribution.  Among many other things, Bernstein 
proved that $\varepsilon$ satisfies
\begin{equation}
\varepsilon > {1\over 4(1+\sqrt 2\kern 1pt)(2\kern .3pt n-1)} 
\end{equation}
for any $n\ge 1$, which implies, since $4(1+\sqrt 2\kern 1pt) \approx 9.66$,
\begin{equation}
n > {1\over 20\kern .5pt \varepsilon}.
\label{bern}
\end{equation}
Since $\varepsilon < 1/2$ in (\ref{prob2}) implies $n\ge 1$, this
establishes condition (\ref{Nbound}) of Theorem~\ref{thm2}.

To establish condition (\ref{cbound}), we make use of (\ref{bern}).
Given $\varepsilon$, let 
$n$ and $\{c_k^{}\}$ define a solution (\ref{prob}) of Problem~E.\ \ If we
split the series into roughly the first quarter and the last three-quarters,
\begin{equation}
\sum_{k=0}^{\lfloor 1/(80\kern .3pt \varepsilon)\rfloor} c_k^{} x^{2k} +
\sum_{k=1 + \lfloor 1/(80\kern .3pt \varepsilon)\rfloor}^n c_k^{} x^{2k} ,
\label{parts}
\end{equation}
then by (\ref{bern}), the first part can approximate $|x|$ no more
closely than $4\kern .5pt \varepsilon$.  More to our purpose, by a linear scaling,
it can approximate $|x|$ over the subinterval $[-1/2,1/2]$ no more
closely than $2\kern .5pt \varepsilon$.  Therefore, since the sum of
the two series in (\ref{parts}) has accuracy better than $\varepsilon$,
the second series must have
maximal size at least $\varepsilon$ over $[-1/2,1/2]$.  Since $|x^{2k}| \le 2^{-2k}$
for $x\in [-1/2,1/2]$, this implies that there must be some
huge coefficients.  Specifically, summing a power
series involving powers of $4$ shows that the second series of (\ref{parts}) is bounded by
\begin{equation}
\sum_{k=1 + \lfloor 1/(80\kern .3pt \varepsilon)\rfloor}^n | c_k^{} x^{2k}| \le
{4\over 3} \kern1pt \cmax\kern 1pt  2^{-2 (\lfloor 1/(80\kern .3pt \varepsilon)\rfloor + 1)}.
\label{pbound}
\end{equation}
Therefore we must have
\begin{equation}
{4\over 3} \kern1pt \cmax\kern 1pt  2^{-2 (\lfloor 1/(80\kern .3pt \varepsilon)\rfloor + 1)} > \varepsilon,
\end{equation}
that is, 
\begin{equation}
\cmax > 0.75\kern .7pt \varepsilon \kern 0.8pt 2^{2 (\lfloor 1/(80\kern .3pt \varepsilon)\rfloor + 1)}.
\end{equation}
This implies (\ref{cbound}), completing the proof of Theorem~\ref{thm2}. \qed

\section{Numerical estimates}

The theorem and proof just given were all about lower bounds, but now let us look
at more accurate (though unrigorous) estimates.  Bernstein~\cite{bernstein} also proved
that the best degree $\twon$ maximum-norm
approximation errors $\varepsilon$ satisfy
\begin{equation}
\varepsilon \sim {\beta\over \twon}, \quad n\to\infty
\end{equation}
for some $\beta$, and in 1985 Varga and Carpenter~\cite{vc} gave the numerical estimate
\begin{equation}
\beta \approx 0.28016949902386913303643649\dots.
\label{estimate}
\end{equation}
To achieve
$\varepsilon \le 10^{-1}$, $10^{-2}$, $10^{-3}$, and $10^{-4}$, respectively,
this suggests (rounding up to the
next even numbers) that we will
need degrees $\twon$ of approximately $4$, $28$, $282$, and $2802$.
It turns out that the actual minimal degrees (as computed with the
Chebfun {\tt minimax} command~\cite{chebfun,minimax}) are exactly these four numbers. 
For accuracy $10^{-6}$, for example, though this is beyond Chebfun, it seems clear that
the required degree will be close to $n=280{,}170$.

Thus we see again that an approximation (\ref{prob}) requires
degrees of order $O(1/\varepsilon)$, but why are the coefficients so large?
The explanation is that the
monomials $1,x^2, x^4, \dots , x^{\twon}$ are an exponentially ill-behaved basis for the space
of even degree $\twon$ polynomials on $[-1,1]$.  Numerical analysts quantify this observation
by noting that the condition number of this set of functions is of the
approximate order
\begin{equation}
\kappa_{\twon} \approx (1+\sqrt 2 \kern 1pt)^{\twon} \approx 10^{\kern 0.5pt 0.766 n}
\end{equation}
\cite{beck,gautschi}.
With $\twon =280{,}170$ for accuracy $10^{-6}$,
this suggests the expansion coefficients will need to be of order
about $10^{107{,}000}$.  
Our best empirical approximation based on
calculations for $n$ up to $300$ is
\begin{equation}
\max_k^{} |c_k^{}| \approx 0.066\times {(1+\sqrt 2\kern 1pt)^{\twon}\over n^{1.5}}.
\label{model}
\end{equation}
Table 1 summarizes our computations and estimates for
accuracies $\varepsilon = 10^{-1}, \dots, 10^{-8}.$

\begin{table}[h]
\begin{center}
\caption{Numerical data for approximation of $x$ by even powers on $\intt$, or
equivalently, approximation of $|x|$ on $[-1,1]$.  The numbers
up to $2802$ are based on numerical computations, and the remaining ones
are estimates.}\label{table}%
\begin{tabular}{@{}rrl@{}}
\toprule
accuracy $\varepsilon$ & minimal degree $\twon$ & maximal coefficient $|c_k^{}|$ \\
\midrule
$10^{-1}~~~$ &        $4~~~~~$ & $~~~~~~~~1.93$ \\
$10^{-2}~~~$ &       $28~~~~~$ & $~~~~~~~~7.4\times 10^{7}$\\
$10^{-3}~~~$ &      $282~~~~~$ & $~~~~~~~~3.5\times 10^{103}$\\
$10^{-4}~~~$ &     $2802~~~~~$ & $~~~~~~~~10^{1068}$ \\
$10^{-5}~~~$ &    $28018~~~~~$ & $~~~~~~~~10^{10{,}700}$ \\
$10^{-6}~~~$ &   $280170~~~~~$ & $~~~~~~~~10^{107{,}000}$ \\
$10^{-7}~~~$ &  $2801696~~~~~$ & $~~~~~~~~10^{1{,}070{,}000}$ \\
$10^{-8}~~~$ & $28016950~~~~~$ & $~~~~~~~~10^{10{,}700{,}000}$ \\
\botrule
\end{tabular}
\end{center}
\end{table}

\begin{figure}[h]
\vspace{8pt}
\begin{center}
\includegraphics[scale=.65]{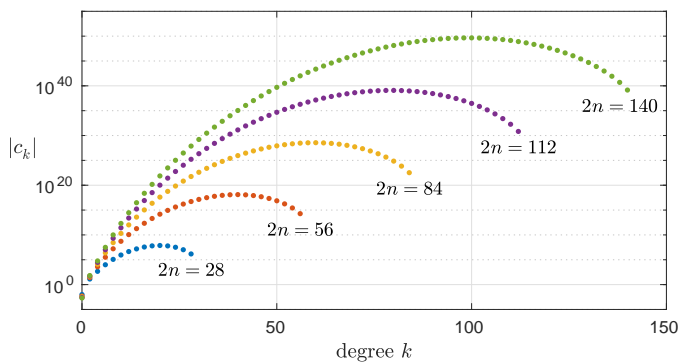}
\caption{Monomial coefficients $|c_k^{}|$ of best approximations of $x$ by even powers on $\intt$, or
equivalently, of $|x|$ on $[-1,1]$, for $\twon = 28, 56, \dots, 140$.  These values correspond
to approximation errors approximately $0.01, 0.01/2, \dots, 0.01/5$.  The end of
the curve is of order $2^{2n}$, and the peak in the middle is of
order $(1+\sqrt 2\kern 1pt)^{2n}$.\label{fig1}}
\end{center}
\end{figure}

Figure~\ref{fig1} illustrates graphically where the big coefficients
come from.  For $\twon = 28, 56, \dots, 140$, it plots the coefficients
$|c_k^{}|$ for $k = 0, 1, \dots, n$ in a monomial expansion of the
best approximations.

\section{A remark about mathematics}

Theorem~\ref{thm2} is startling and interesting.  From the usual
mathematical point of view, however, it is
not much more than that.  After all, M\"untz's theorem remains
valid and beautiful.  From the usual mathematical perspective,
M\"untz's theorem expresses a fundamental truth, and Theorem~\ref{thm2},
however interesting, is an engineering footnote.

As I have discussed in the context of other problems~\cite{yogi,apology},
I believe this usual perspective is too comfortable.
Theorem~\ref{thm2} implies that typical 
sets of powers deemed useful by M\"untz's theorem would
in fact be useless in any actual application.  If it was not the business
of mathematicians to notice and analyze this effect in the past 110 years,
then whose business was it?

\backmatter

\bmhead{Acknowledgments}

I am grateful for helpful suggestions to Michael
Ganzburg, Daan Huybrechs, Christian Lubich, Doron Lubinsky, Yuji Nakatsukasa,
Allan Pinkus, and Endre S\"uli.

\bmhead{Conflict of Interest}

The author has no relevant financial or non-financial interests to disclose.

\end{document}